\documentclass[fleqn,11pt,a4paper,leqno]{article}

\usepackage{graphicx}
\usepackage{graphics}
\usepackage{latexsym}
\usepackage{amssymb,amsmath}
\usepackage{multicol}
\usepackage{bussproofs}
\usepackage{enumerate}
\usepackage{xcolor}
\usepackage[latin1]{inputenc}
\usepackage{verbatim, color}
\usepackage{enumerate}

\title{proof of associative of $\land$}

\author{}

\date{}

\numberwithin{equation}{section}

\newtheorem{definicion}{Definition}[section]

\newtheorem{definition}[definicion]{Definition}

\newtheorem{Theorem}[definicion]{Theorem}

\newtheorem{Corollary}[definicion]{Corollary}

\newtheorem{lema}[definicion]{Lemma}

\newenvironment{Proof}{\noindent\bf Proof \rm}{$\hfill
\square$}

\begin{document}

%%MATT's  "thlist" environment
 \newcounter{thlistctr}
 \newenvironment{thlist}{\
 \begin{list}%
 {\alph{thlistctr}}%
 {\setlength{\labelwidth}{2ex}%
 \setlength{\labelsep}{1ex}%
 \setlength{\leftmargin}{6ex}%
 \renewcommand{\makelabel}[1]{\makebox[\labelwidth][r]{\rm (##1)}}%
 \usecounter{thlistctr}}}%
 {\end{list}}
%%END MATT's "thlist"

	\title{Semidistributivity and Whitman Property\\
	 in\\ Implication Zroupoids }
	
	\author{Juan M. CORNEJO and Hanamantagouda P. SANKAPPANAVAR}

\maketitle

\begin{abstract}

%The main result of this paper is the every implication zroupoid satisfies semi-distributive identity.

%In \cite{sankappanavarMorgan2012}, the second-named author defined and initiated the investigations into %studied in \cite{sankappanavarMorgan2012} t
%the variety $\mathcal{I}$ of algebras, called implication zroupoids %($\mathcal{I}$-zroupoids) 
%that generalize De Morgan algebras.  
%\begin{comment}
%In 1934, Bernstein gave a system of axioms for Boolean algebras in terms of implication only. %however, 
% In 2012, the second author of this paper extended Bernstein's theorem
%to De Morgan algebras in \cite{sankappanavarMorgan2012} by  
% showing %in \cite{sankappanavarMorgan2012}
%  that the varieties of De Morgan algebras, Kleene algebras, and  Boolean algebras are term-equivalent,  
%to varieties whose defining axioms use only the implication $\to$ and the constant $0$. 
   
%The primary role played by the identity (I): $(x \to y) \to z \approx [(z' \to x) \to (y \to z)']'$, where $x' : = x \to 0$, which
%occurs as an axiom in the definition of each of those new varieties 
%motivated the second author of this paper to 
%introduce a new (equational) class of algebras called ``Implication zroupoids'' in   
%\cite{sankappanavarMorgan2012}.
%\end{comment}

In 2012, the second author introduced and studied in \cite{sankappanavarMorgan2012} the variety $\mathcal{I}$ of implication zroupoids 
that generalize De Morgan algebras and $\lor$-semilattices with $0$.
An algebra $\mathbf A = \langle A, \to, 0 \rangle$, where $\to$ is binary and $0$ is a constant, is called an 
\emph{implication zroupoid} ($\mathcal{I}$-zroupoid, for short) if $\mathbf A$ satisfies:
 $(x \to y) \to z \approx [(z' \to x) \to (y \to z)']'$, where $x' : = x \to 0$,  and  
 $ 0'' \approx 0$.   Let $\mathcal{I}$ denote the variety of implication zroupoids   
%\begin{comment}
% For details on the motivation leading to these algebras, we refer the reader to \cite{sankappanavarMorgan2012} (or the relevant papers mentioned at the end of this paper.    
%The investigations into the structure of the lattice of subvarieties of $\mathcal{I}$, begun in \cite{sankappanavarMorgan2012}, have continued in \cite{cornejo2016order, cornejo2016semisimple, cornejo2015implication}. %to the investigation of the structure of the lattice of subvarieties of $\mathcal{I}$, and to making further contributions to the theory of implication zroupoids.
%This paper is a sequel to the above-mentioned papers.   
%\end{comment}
and $\mathbf A \in \mathcal{I}$.  For $x,y \in \mathbf A$, let $x \land y := (x \to y')'$ and  $x \lor y := (x' \land y')'$.  In an earlier paper we had proved that if $\mathbf A \in \mathcal{I}$, then the algebra $\mathbf A_{mj} = \langle A, \lor, \land \rangle$ is a bisemigroup.
%The identity (BR): $x \land (x \lor y) \approx x \lor (x \land y)$ is called the Birkhoff's identity.  
In this paper we generalize the notion of semi-distributivity from lattices to bisemigroups and
%The main purpose of this paper is to 
prove that, for every $\mathbf A \in \mathcal{I}$, the bisemigroup  
$\mathbf A_{mj}$ is semidistributive. % {\color{blue}(To be modified.)}
%the variety $\mathcal{I}$ satisfies the Birkhoff's identity, thus extending an earlier result of the authors, proved in \cite{cornejo2016derived}.  It follows from this result that there are bisemigroups that are not bisemilattices but that satisfy the Birkhoff's identity, which suggests a more general notion, than Birkhoff systems, of ``Birkhoff bisemigroups'' as bisemigroups satisfying the Birkhoff identity. 
Secondly, we generalize the Whitman Property from lattices to bisemigroups and prove that the subvariety $\mathcal{MEJ}$ of $\mathcal I$, defined by the identity: $x \land y \approx x \lor y$, satisfies the Whitman Property.

\end{abstract}
%{\color{blue}
\section{Introduction}
% To be modified}
\
%{\color{blue}

 Bernstein \cite{bernstein1934postulates} gave a system of axioms, in 1934, for Boolean algebras in terms of implication only. %however, his original axioms were not equational.  %A quick look at his axioms would reveal that, 
% with an additional constant, they 
%could easily be translated into equational ones. 
The second author of this paper extended Bernstein's theorem
to De Morgan algebras in \cite{sankappanavarMorgan2012} by  
 showing %in \cite{sankappanavarMorgan2012}
  that the varieties of De Morgan algebras, Kleene algebras, and  Boolean algebras are term-equivalent,  
to varieties whose defining axioms use only the implication $\to$ and the constant $0$.    
The primary role played by the identity (I): $(x \to y) \to z \approx [(z' \to x) \to (y \to z)']'$, where $x' : = x \to 0$, which
occurs as an axiom in the definition of each of those new varieties 
led him, in 2012, to 
introduce a new (equational) class of algebras called ``Implication zroupoids'' in   
\cite{sankappanavarMorgan2012}.

An algebra $\mathbf A = \langle A, \to, 0 \rangle$, where $\to$ is binary and $0$ is a constant, is called an 
\emph{implication zroupoid} ($\mathcal{I}$-zroupoid, for short) if $\mathbf A$ satisfies:

\begin{itemize}
\item[{\rm(I)}]  $(x \to y) \to z \approx [(z' \to x) \to (y \to z)']'$, where $x' : = x \to 0$,  and  
\item[{\rm(I$_0$)}]  $ 0'' \approx 0$.  
\end{itemize}
 Let $\mathcal{I}$ denote the variety of implication zroupoids.  These algebras generalize De Morgan algebras and $\lor$-semilattices with zero.  For more details on the motivation leading to these algebras, we refer the reader to \cite{sankappanavarMorgan2012} (or the relevant papers mentioned at the end of this paper).  
 
% For detailed summary of the results contained in the above-mentioned papers, we refer the reader to the introduction of \cite{cornejo2016derived}.
 
The investigations into the (complex) structure of the lattice of subvarieties of $\mathcal{I}$, begun in \cite{sankappanavarMorgan2012}, have continued in \cite{cornejo2016order}, \cite{cornejo2016semisimple}, \cite{cornejo2015implication}, \cite{cornejo2016derived}, \cite{cornejo2016Bol-Moufang}, \cite{associativetype},  
\cite{cornejo2016weakassociative}, \cite{Implsem} and \cite{cornejo2020birkhoff}.  
The present paper is a sequel to this series of papers and is devoted to %investigating the structure of the lattice of subvarieties of $\C{I}$, and also to 
making further contributions to the theory of implication zroupoids. %In this paper we are interested in the 

Throughout this paper we use the following definitions:
\[
	\textup{(M)} \quad 	x \land y := (x \to y')' \quad\text{ and }\quad 
	\textup{(J)} \quad  x \lor y := (x' \land y')'.
\]
%{\color{red} We will denote $x \leq y$ if $x \wedge y = x$.}

The relation $\leq$ is defined as follows: $x \leq y$ if and only if $x \wedge y = x$.  We note that, in general, it is not a partial order on $\mathcal{I}$.  However, $\mathcal{I}_{2,0}$ is the maximal subvariety of $\mathcal{I}$ in which the relation $\leq$ is a partial order (see \cite{cornejo2016order}). \\

With each  $\mathbf A \in \mathcal{I}$, we associate the following algebra:
\[
	\mathbf{A{^{mj}}} := \langle A, \land, \lor, 0 \rangle\\	
	\mathbf{A{_{mj}}} := \langle A, \land, \lor \rangle.   
\]

It was proved in \cite[Corollary 4.6]{cornejo2016derived} that if $\mathbf{A}$ is an implication zroupoid, then 
$\mathbf{A{_{mj}}}$ is a bisemigroup (i.e., an algebra with two binary operations which are both associative.)

\begin{Theorem} {\rm \cite[Corollary 4.6]{cornejo2016derived}}  \label{Theo_Assoc}
	If $\mathbf A \in \mathcal I$ then $\langle A, \wedge \rangle$ and $\langle A, \vee \rangle$ are semigroups.
\end{Theorem}

Two of the important subvarieties of $\mathcal{I}$ are:
 $\mathcal{I}_{2,0}$ and $\mathcal{MC}$ which are defined relative to $\mathcal{I}$, respectively, by the following identities, where $x \land y := (x \to y')' $ :      
\begin{equation} \label{eq_I20}  \tag{{\rm I}$_{2,0}$}
x'' \approx x,
\end{equation}
\begin{equation} \label{eq_MC} \tag{MC}
x \wedge y \approx y \wedge x.
\end{equation}

%\begin{definition}
	%Let $\mathbf{A} \in \mathcal{I}$.  
%We say $\mathbf{A}$ is {\rm involutive} if $A \in 
Members of the variety $\mathcal{I}_{2,0}$ are called \emph {involutive}, and members of   %$\mathbf{A}$ is {\rm meet-commutative} if $A \in 
$\mathcal{MC}$ are called \emph {meet-commutative}.  An algebra $\mathbf{A} \in  \mathcal{I}$ is {\rm symmetric} if $\mathbf{A}$ is both involutive and meet-commutative.    Let $\mathcal{S}$ denote the variety of symmetric $\mathcal{I}$-zroupoids.  In other words, $\mathcal{S} = \mathcal{I}_{2,0} \cap \mathcal{MC}$.

%{\color{red} Sanka, I think that we need to work in the introduction a bit more.}

%\end{definition}
%Another important subvariety of $\mathcal{I}$ is the variety $\mathbb{S}$ of $\mathcal{I}$ of symmetric implication zroupoids. 
%Recall $\mathbb{S}:= \mathbb{I}_{2,0} \cap \mathbb{MC}$.  Members of $\mathbb{S}$ are called symmetric zroupoid. 

The notions of meet-semidistributivity, join-semidistributivity and semidistributivity for lattices were first defined, in 1961, by J\'{o}nsson \cite{Jo61}, who proved that free lattices are semidistributive.  These notions have been investigated in group theory and in semigroup theory also.  For example, It was shown by Shiryaev \cite{Sh85} 
that meet-semidistributivity and distributivity are equivalent for the lattice of subgroups of a group; on the other hand, join-semidistributivity and distributivity are distinct. 
In 1999, \cite{JoJo99} showed that each of meet-semidistibutivity and join-semidistributivity is equivalent to distributivity on the lattice of inverse subsemigroups containing all the idempotents of an inverse semigroup.

It is clear that the notion of semidistributivity extends naturally from lattices to bisemigroups.

	\begin{definition} A bisemigroup is
		\begin{itemize}
			\item meet-semidistributive if it satisfies the following conditions:\\
		{\rm($M_1$)}	\quad $x \land y \approx x \land z$ implies $x \land (y \lor z) \approx x \land y$ (left meet-semidistributive law), \\ 
		{\rm($M_2$)}	\quad $ x \land y \approx z \land y$ implies  $(x \lor z) \land y \approx x \land y$ (right meet-semidistributive law)
			\item join-semidistributive if it satisfies  the following conditions:\\
			{\rm($J_1$)}	\quad  $x \lor y \approx x \lor z$ implies $x \lor (y \land z) \approx x \lor y$  (left join-semidistributive law),\\
			{\rm($J_2$)}   \quad $x \lor y \approx z \lor y$ implies $(x \land z) \lor y \approx x \lor y$  (right join-semidistributive law).
			\item semidistributive if it is both meet-semi-distributive and join-semi-distributive.
		\end{itemize}
	\end{definition}

%{\color{red} HABLAR DE LA DEFINICION! LOCALIDAD DE LAS IDENTIDADES???}

It is easy to see that semilattices, viewed as bisemilattices, where the two binary operations are the same, are semidistributive; and De Morgan algebras are clearly semidistributive.  The variety $\mathcal{I}$ of implication zroupoids contains both the varieties of $\lor$-semilattices with $0$ and De Morgan algebras.  Moreover, every implication zroupoid $\mathbf A$ gives rise to the bisemigroup $\mathbf A_{mj}$ .  So, it is natural to ask whether $\mathbf A_{mj}$ is semidistributive, for an I-zroupoid $\mathbf A$.  %The main result answers this question in the positive.}
The purpose of this note is to answer this question in the positive by proving that $A_{mj}$ is semidistributive for every I-zroupoid $\mathbf A$. 

\section{Preliminaries}
\
We refer the reader to the books \cite{balbesDistributive1974}, \cite{burrisCourse1981} and \cite{Ra74} for the concepts and results assumed in this paper.

We now present some preliminary results that will be useful later. 
\begin{lema} {\rm \cite[Theorem 8.15]{sankappanavarMorgan2012}} \label{general_properties_equiv}
	%Let $\mathbf A$ be a I-zroupoid. Then 
	The following identities are equivalent in the variety $\mathcal I$:
	\begin{enumerate}[{\rm (a)}]
		\item $0' \to x \approx x$, \label{TXX} %True implies x = x
		\item $x'' \approx x$,
		\item $(x \to x')' \approx x$, \label{reflexivity}
		\item $x' \to x \approx x$. \label{LeftImplicationwithtilde}
	\end{enumerate}
\end{lema}

\begin{lema} \label{lemaOfIZI} {\bf \cite[Lemma 3.4]{cornejo2015implication}}  
	Let $\mathbf A$ be an I-zroupoid.  Then $\mathbf A$ satisfies:
	\begin{enumerate}[{\rm (a)}] \label{Lemma_300315_01}
		\item $(x \to y) \to z \approx [(x \to y) \to z]'' $, \label{lemaOfIZI3_4}  \label{060415_01}
		\item $(x \to y)' \approx (x'' \to y)' $. \label{lemaOfIZI2_4} \label{300315_02}
	\end{enumerate}
\end{lema}

\begin{lema}\label{general_properties} \label{general_properties2}
	Let $\mathbf A \in \mathcal I_{2,0}$. Then $\mathbf A$ satisfies:
	\begin{enumerate}[{\rm (1)}]
		\item $(x \to 0') \to y \approx (x \to y') \to y$, \label{281014_05}
%		\item $x' \to 0' \approx 0 \to x$, \label{cuasiConmutativeOfImplic2}
		\item $0 \to x' \approx x \to 0'$, \label{cuasiConmutativeOfImplic}
%		\item $(x \to 0') \to (y \to z) \approx ((0 \to x) \to y) \to z$, \label{260615_01}
%		\item $(0 \to x) \to (0 \to y) \approx x \to (0 \to y)$, \label{311014_06}
		\item $0 \to (x \to y) \approx x \to (0 \to y)$, \label{071114_04}
%		\item $0 \to (x' \to y)' \approx x \to (0 \to y')$, \label{071114_01}
%		\item $(x \to y) \to (y \to z) \approx (0 \to x') \to (y \to z)$, \label{250615_04}
%		\item $((x \to y) \to z) \to (z \to u) \approx (0 \to x) \to ((y \to z) \to (z \to u))$, \label{250615_06} 
%		\item $x \to y \approx x \to (x \to y)$, \label{031114_04} %(246)
		\item $(y \to x) \to y \approx (0 \to x) \to y$, \label{291014_10}
%		\item $(x \to y)' \to (0 \to x)' \approx y' \to x'$, \label{271114_02} %{\color{red} IDENTIDAD 557}
%		\item $(x \to y)' \to y \approx x \to y$, \label{250615_01} 
%		\item $[x \to (y \to x)']' \approx (x \to y) \to x$, \label{291014_09}
%		\item $x \to ((0 \to x) \to y) \approx x \to y$, \label{250615_22} 
%		\item $x \to (y \to x') \approx y \to x'$, \label{281114_01}
%		\item $(x \to y) \to y' \approx y \to (x \to y)'$, \label{071114_05}
%		\item $((x \to y) \to (z \to x)) \to u \approx (y \to 0') \to ((z \to x) \to u)$, \label{260615_05} 
%		\item $(z \to x) \to (y \to z) \approx (0 \to x) \to (y \to z)$, \label{080415_01}
		\item $(x \to y')' \to z \approx x \to (y \to z)$, \label{140715_20}
		\item $0 \to (x \to y')' \approx 0 \to (x' \to y)$, \label{191114_05}
%		\item $0 \to (0 \to x)' \approx 0 \to x'$, \label{031114_07}
%		\item $[x' \to (0 \to y)]' \approx (0 \to x) \to (0 \to y)'$. \label{031114_06}
\item $[x \to (x' \to y)']' \approx x' \to (0 \to y')'$.  \label{130919_04}
	\end{enumerate}
\end{lema}

\begin{Proof}
	Item (\ref{281014_05}) can be found in \cite[Lemma 2.7]{cornejo2016derived}.
	Item
	%s (\ref{cuasiConmutativeOfImplic2}) and 
	(\ref{cuasiConmutativeOfImplic}) is proved in \cite{sankappanavarMorgan2012}.
	The proofs of items 
	%(\ref{311014_06}), 
	(\ref{071114_04}),   
	%(\ref{071114_01}), (\ref{031114_04}), 
	(\ref{291014_10}), and
	%(\ref{271114_02}), (\ref{291014_09}), (\ref{281114_01}), (\ref{071114_05}), 
	(\ref{191114_05}) 
	%(\ref{031114_07}), (\ref{031114_06}) 
	can be found in \cite{cornejo2016order}.
	Items 
	%(\ref{260615_01}), (\ref{250615_04}), (\ref{250615_06}), (\ref{250615_01}), (\ref{250615_22}), (\ref{260615_05}), (\ref{080415_01}), 
	(\ref{140715_20}) and
	%The item 
	(\ref{130919_04}) are proved, respectively, in  \cite{cornejo2016derived} and \cite{cornejo2020birkhoff}.
\end{Proof}

%\vspace{1cm}
\section{Semidistributivity of involutive I-zroupoids}
\
%The notion of semidistributivity was first defined, in 1961, for lattices by J\'{o}nsson \cite{Jo61}, who  proved that free lattices are semi-distributive.  We find it necessary to extend this notion to bisemigroups.
%\begin{definition}
%\noindent A bisemigroup is meet-semi-distributive if it satisfies:\\
%(SD$_\land$) \quad $x \land y \approx x \land z$ implies $x \land (y \lor z) \approx x \land y$.  Dually, it is join-semi-distributive if it satisfies:\\
%(SD$_\lor$)   \quad $x \lor y \approx x \lor z$ implies $x \lor (y \land z) \approx x \lor y$.\\
%It is semi-distributive if it is both meet-semi-distributive and join-semi-distributive. 
%An I-zroupoid is meet-semi-distributive (join-semi-distributive, semi-distributive, respectively) if $\mathbf A_{mj}$ is.
%\end{definition}

\indent  In this section we prove that if  $\mathbf{A} \in \mathcal I_{2,0}$, then the bisemigroup $\mathbf{A}_{mj}$ is semi-distributive.   This result will play an important role in the proof of the main Theorem (Theorem \ref{Main}) in the next section.  To this end, we need the following lemmas.

\begin{lema} \label{Lemma_070220_01}
	Let $\mathbf A \in \mathcal I_{2,0}$ and $a,b \in A$ such that $a \to b' = a \to c'$. Then
	\begin{enumerate}[{\rm (1)}]
		\item $(a' \to b) \to a' = (a' \to c) \to a'$, \label{070220_01}
		\item $0 \to (a' \to b)' = 0 \to (a' \to c)'$, \label{070220_02}
		\item $a' \to (b \to 0')' = a' \to (c \to 0')'$, \label{070220_03}
		\item $(b' \to c) \to a' = c \to a'$, \label{070220_04}
		\item $a \to (b' \to c)' = a \to b'$. \label{070220_05}
	\end{enumerate}
\end{lema}

\begin{Proof}
\begin{enumerate}
	\item[(\ref{070220_01})]  
%$$
%\begin{array}{lcll}
%(a' \to b) \to a' & = & (0 \to b) \to a' & \mbox{by Lemma \ref{general_properties2} (\ref{291014_10})} \\
%& = & (b' \to 0') \to a' & \mbox{by Lemma \ref{general_properties} (\ref{cuasiConmutativeOfImplic})  and $x'' \approx x$}\\
%& = & [(a'' \to b') \to (0' \to a')']' & \mbox{by (I)} \\
%& = & [(a \to b') \to (0' \to a')']' & \mbox{} \\
%& = & [(a \to c') \to (0' \to a')']' & \mbox{by hypothesis} \\
%& = & [(a'' \to c') \to (0' \to a')']' & \mbox{} \\
%& = & (c' \to 0') \to a' & \mbox{by (I)} \\
%& = & (0 \to c) \to a' & \mbox{by Lemma \ref{general_properties} (\ref{cuasiConmutativeOfImplic}) and $x'' \approx x$} \\
%& = & (a' \to c) \to a' & \mbox{by Lemma \ref{general_properties2} (\ref{291014_10})} 
%\end{array}
%$$
\noindent $(a' \to b) \to a' $
$\overset{   \ref{general_properties2} (\ref{291014_10}) 
}{=}  (0 \to b) \to a' $
$\overset{   \ref{general_properties} (\ref{cuasiConmutativeOfImplic})  \ {\rm and}\  x'' \approx x
}{=}  (b' \to 0') \to a' $
$\overset{  (I) 
}{=}  [(a'' \to b') \to (0' \to a')']' $
$\overset{  
}{=}  [(a \to b') \to (0' \to a')']' $
$\overset{ hyp 
}{=}  [(a \to c') \to (0' \to a')']' $
$\overset{  
}{=}  [(a'' \to c') \to (0' \to a')']' $
$\overset{  (I) 
}{=}  (c' \to 0') \to a' $
$\overset{   \ref{general_properties} (\ref{cuasiConmutativeOfImplic}) \ {\rm and}\   x'' \approx x 
}{=}  (0 \to c) \to a' $
$\overset{   \ref{general_properties2} (\ref{291014_10}) 
}{=}  (a' \to c) \to a' $.

	\item[(\ref{070220_02})]  
%$$
%\begin{array}{lcll}
%0 \to (a' \to b)' & = & 0 \to (a'' \to b') & \mbox{by Lemma \ref{general_properties2} (\ref{191114_05}) and $x'' \approx x$} \\
%& = & 0 \to (a \to b') & \mbox{} \\
%& = & 0 \to (a \to c') & \mbox{by hypothesis} \\
%& = & 0 \to (a'' \to c') & \mbox{} \\
%& = & 0 \to (a' \to c)' & \mbox{by Lemma \ref{general_properties2} (\ref{191114_05})  and $x'' \approx x$} 
%\end{array}
%$$
\noindent $0 \to (a' \to b)' $
$\overset{   \ref{general_properties2} (\ref{191114_05}) \ {\rm and}\  x'' \approx x
}{=}  0 \to (a'' \to b') $
$\overset{  
}{=}  0 \to (a \to b') $
$\overset{ hyp 
}{=}  0 \to (a \to c') $
$\overset{  
}{=}  0 \to (a'' \to c') $
$\overset{   \ref{general_properties2} (\ref{191114_05})  \ {\rm and}\  x'' \approx x 
}{=}  0 \to (a' \to c)' $.

	\item[(\ref{070220_03})] 
%$$
%\begin{array}{lcll}
%a' \to (b \to 0')' & = & a' \to (0 \to b')' & \mbox{by Lemma \ref{general_properties} (\ref{cuasiConmutativeOfImplic})} \\
%& = & [a \to (a' \to b)']' & \mbox{by Lemma \ref{general_properties} (\ref{130919_04})} \\
%& = & [a'' \to (a' \to b)']' & \mbox{} \\
%& = & [(a' \to 0) \to (a' \to b)']' & \mbox{} \\
%& = & [((a' \to b)'' \to a') \to (0 \to (a' \to b)')']'' & \mbox{by (I)} \\
%& = & ((a' \to b) \to a') \to (0 \to (a' \to b)')' & \mbox{} \\
%& = & ((a' \to c) \to a') \to (0 \to (a' \to b)')' & \mbox{by (\ref{070220_01})} \\
%& = & ((a' \to c) \to a') \to (0 \to (a' \to c)')' & \mbox{by (\ref{070220_02})} \\
%& = & [((a' \to c)'' \to a') \to (0 \to (a' \to c)')']'' & \mbox{} \\
%& = & [(a' \to 0) \to (a' \to c)']' & \mbox{by (I)} \\
%& = & [a \to (a' \to c)']' & \mbox{} \\
%& = & a' \to (0 \to c')' & \mbox{by Lemma \ref{general_properties} (\ref{130919_04})} \\
%& = & a' \to (c \to 0')' & \mbox{by Lemma \ref{general_properties} (\ref{cuasiConmutativeOfImplic})} 
%\end{array}
%$$
\noindent $a' \to (b \to 0')' $
$\overset{   \ref{general_properties} (\ref{cuasiConmutativeOfImplic}) 
}{=}  a' \to (0 \to b')' $
$\overset{   \ref{general_properties} (\ref{130919_04}) 
}{=}  [a \to (a' \to b)']' $
$\overset{  
}{=}  [a'' \to (a' \to b)']' $
$\overset{  
}{=}  [(a' \to 0) \to (a' \to b)']' $
$\overset{  (I) 
}{=}  [((a' \to b)'' \to a') \to (0 \to (a' \to b)')']'' $
$\overset{  
}{=}  ((a' \to b) \to a') \to (0 \to (a' \to b)')' $
$\overset{  (\ref{070220_01}) 
}{=}  ((a' \to c) \to a') \to (0 \to (a' \to b)')' $
$\overset{  (\ref{070220_02}) 
}{=}  ((a' \to c) \to a') \to (0 \to (a' \to c)')' $
$\overset{  
}{=}  [((a' \to c)'' \to a') \to (0 \to (a' \to c)')']'' $
$\overset{  (I) 
}{=}  [(a' \to 0) \to (a' \to c)']' $
$\overset{  
}{=}  [a \to (a' \to c)']' $
$\overset{   \ref{general_properties} (\ref{130919_04}) 
}{=}  a' \to (0 \to c')' $
$\overset{   \ref{general_properties} (\ref{cuasiConmutativeOfImplic}) 
}{=}  a' \to (c \to 0')' $.

	\item[(\ref{070220_04})] 
%$$
%\begin{array}{lcll}
%(b' \to c) \to a' & = & [(a'' \to b') \to (c \to a')']' & \mbox{by (I)} \\
%& = & [(a \to b') \to (c \to a')']' & \mbox{} \\
%& = & [(a \to c') \to (c \to a')']' & \mbox{by hypothesis} \\
%& = & [(a'' \to c') \to (c \to a')']' & \mbox{} \\
%& = & (c' \to c) \to a' & \mbox{by (I)} \\
%& = & c \to a' & \mbox{by Lemma \ref{general_properties_equiv} (\ref{LeftImplicationwithtilde})} 
%\end{array}
%$$
\noindent $(b' \to c) \to a' $
$\overset{  (I) 
}{=}  [(a'' \to b') \to (c \to a')']' $
$\overset{  
}{=}  [(a \to b') \to (c \to a')']' $
$\overset{ hyp 
}{=}  [(a \to c') \to (c \to a')']' $
$\overset{  
}{=}  [(a'' \to c') \to (c \to a')']' $
$\overset{  (I) 
}{=}  (c' \to c) \to a' $
$\overset{   \ref{general_properties_equiv} (\ref{LeftImplicationwithtilde}) 
}{=}  c \to a' $.

	\item[(\ref{070220_05})] 
%$$
%\begin{array}{lcll}
%a \to (b' \to c)' & = & a'' \to (b' \to c)'  & \mbox{} \\
%& = & (a' \to 0) \to (b' \to c)' & \mbox{} \\
%& = & [((b' \to c)'' \to a') \to (0 \to (b' \to c)')']' & \mbox{by (I)} \\
%& = & [((b' \to c) \to a') \to (0 \to (b' \to c)')']' & \mbox{} \\
%& = & [(c \to a') \to (0 \to (b' \to c)')']' & \mbox{by (\ref{070220_04})} \\
%& = & [(c \to a') \to (0 \to (b \to c'))']' & \mbox{by Lemma \ref{general_properties2} (\ref{191114_05})} \\
%& = & [(c \to a') \to (b \to (0 \to c'))']' & \mbox{by Lemma \ref{general_properties2} (\ref{071114_04})} \\
%& = & [(c \to a') \to ((b \to 0')' \to c')']' & \mbox{by Lemma \ref{general_properties2} (\ref{140715_20})} \\
%& = & [(c'' \to a') \to ((b \to 0')' \to c')']' & \mbox{} \\
%& = & (a' \to (b \to 0')') \to c' & \mbox{by (I)} \\
%& = & (a' \to (c \to 0')') \to c' & \mbox{by (\ref{070220_03})} \\ 
%& = & [(c'' \to a') \to ((c \to 0')' \to c')']' & \mbox{by (I)} \\
%& = & [(c'' \to a') \to (c \to (0 \to c'))']' & \mbox{by Lemma \ref{general_properties2} (\ref{140715_20})} \\
%& = & [(c'' \to a') \to (0 \to (c \to c'))']' & \mbox{by Lemma \ref{general_properties2} (\ref{071114_04})} \\
%& = & [(c'' \to a') \to (0 \to (c'' \to c'))']' & \mbox{} \\
%& = & [(c'' \to a') \to (0 \to c')']' & \mbox{by Lemma \ref{general_properties_equiv} (\ref{LeftImplicationwithtilde})} \\
%& = & (a' \to 0) \to c' & \mbox{by (I)} \\
%& = & a \to c' & \mbox{} \\
%& = & a \to b' & \mbox{by hypothesis} 
%\end{array}
%$$
\noindent $a \to (b' \to c)' $
$\overset{  
}{=}  a'' \to (b' \to c)'  $
$\overset{  
}{=}  (a' \to 0) \to (b' \to c)' $
$\overset{  (I) 
}{=}  [((b' \to c)'' \to a') \to (0 \to (b' \to c)')']' $
$\overset{  
}{=}  [((b' \to c) \to a') \to (0 \to (b' \to c)')']' $
$\overset{  (\ref{070220_04}) 
}{=}  [(c \to a') \to (0 \to (b' \to c)')']' $
$\overset{   \ref{general_properties2} (\ref{191114_05}) 
}{=}  [(c \to a') \to (0 \to (b \to c'))']' $
$\overset{   \ref{general_properties2} (\ref{071114_04}) 
}{=}  [(c \to a') \to (b \to (0 \to c'))']' $
$\overset{   \ref{general_properties2} (\ref{140715_20}) 
}{=}  [(c \to a') \to ((b \to 0')' \to c')']' $
$\overset{  
}{=}  [(c'' \to a') \to ((b \to 0')' \to c')']' $
$\overset{  (I) 
}{=}  (a' \to (b \to 0')') \to c' $
$\overset{  (\ref{070220_03})  
}{=}  (a' \to (c \to 0')') \to c' $
$\overset{  (I) 
}{=}  [(c'' \to a') \to ((c \to 0')' \to c')']' $
$\overset{   \ref{general_properties2} (\ref{140715_20}) 
}{=}  [(c'' \to a') \to (c \to (0 \to c'))']' $
$\overset{   \ref{general_properties2} (\ref{071114_04}) 
}{=}  [(c'' \to a') \to (0 \to (c \to c'))']' $
$\overset{  
}{=}  [(c'' \to a') \to (0 \to (c'' \to c'))']' $
$\overset{   \ref{general_properties_equiv} (\ref{LeftImplicationwithtilde}) 
}{=}  [(c'' \to a') \to (0 \to c')']' $
$\overset{  (I) 
}{=}  (a' \to 0) \to c' $
$\overset{  
}{=}  a \to c' $
$\overset{ hyp 
}{=}  a \to b' $.
\end{enumerate}
\end{Proof}

\begin{Theorem}  \label{Theorem_070220_01}
Let $\mathbf A \in \mathcal I_{2,0}$.  Then $\mathbf A$ is meet-semidistributive.
\end{Theorem}

\begin{Proof}  First we will show that $\mathbf A$ satisfies the condition {\rm($M_1$)}. %(SD$_\land$ (left)).}
By hypothesis, we have that $\mathbf A \models x \wedge y \approx x \wedge z$. 
Let $a,b,c \in A$. % such that $a \wedge b = a \wedge c$.  
Now, observe that
%$$
%\begin{array}{lcll}
%a \to b' & = & (a \to b')'' & \mbox{} \\
%& = & (a \wedge b)' & \mbox{} \\
%& = & (a \wedge c)' & \mbox{} \\
%& = & (a \to c')'' & \mbox{} \\
%& = & a \to c' & \mbox{} 
%\end{array}
%$$
\noindent $a \to b' $
$\overset{  
}{=}  (a \to b')'' $
$\overset{  
}{=}  (a \wedge b)' $
$\overset{  
}{=}  (a \wedge c)' $
$\overset{  
}{=}  (a \to c')'' $
$\overset{  
}{=}  a \to c' $.
Hence, by Lemma \ref{Lemma_070220_01} (\ref{070220_05}), we have 
\begin{equation} \label{equation_070220_01}
a \to (b' \to c)' = a \to b'.
\end{equation}
Therefore,
%$$
%\begin{array}{lcll}
%a \wedge (b \vee c)  & = & a \wedge (b' \wedge c')' & \mbox{by definition of } \vee \\
%& = & a \wedge (b' \to c'')'' & \mbox{by definition of }  \wedge  \\
%& = & a \wedge (b' \to c) & \mbox{} \\
%& = & (a \to (b' \to c)')' & \mbox{by definition of } \wedge \\
%& = & (a \to b')' & \mbox{by (\ref{equation_070220_01})} \\
%& = & a \wedge b. & \mbox{by definition of } \wedge 
%\end{array}
%$$
\noindent $a \wedge (b \vee c)  $
$\overset{  \ {\rm def\  of }\  \vee 
}{=}  a \wedge (b' \wedge c')' $
$\overset{  \ {\rm def\  of }\    \wedge  
}{=}  a \wedge (b' \to c'')'' $
$\overset{  
}{=}  a \wedge (b' \to c) $
$\overset{  \ {\rm def\  of }\   \wedge 
}{=}  (a \to (b' \to c)')' $
$\overset{  (\ref{equation_070220_01}) 
}{=}  (a \to b')' $
$\overset{  \ {\rm def\  of }\   \wedge 
}{=}  a \wedge b. $

%{\color{red} 
%Now we desire to check that $\mathbf A$ satisfies (SD$_\land$ (right)). 
%Let us assume then that $\mathbf A \models x \wedge y \approx z \wedge y$ and take $a,b,c \in A$. Hence
%$$
%\begin{array}{lcll}
%(a \vee c) \wedge b & = & (a' \wedge c')' \wedge b & \mbox{by definition of } \vee \\
%& = & ((a' \wedge c')' \to b')' & \mbox{by definition of }  \wedge  \\
%& = & ((a' \to c'')'' \to b')' & \mbox{by definition of }  \wedge  \\
%& = & ((a' \to c) \to b')' & \mbox{} \\
%& = & ((b'' \to a') \to (c \to b')')'' & \mbox{by (I)} \\
%& = & ((b'' \to a') \to (c \wedge b))'' &  \mbox{by definition of }  \wedge  \\
%& = & ((b'' \to a') \to (a \wedge b))'' & \mbox{by hypothesis} \\
%& = & ((b'' \to a') \to (a \to b')')''  & \mbox{by definition of }  \wedge  \\
%& = & ((a' \to a) \to b')' & \mbox{by (I)} \\
%& = & (a \to b')' & \mbox{by Lemma \ref{general_properties_equiv} (\ref{LeftImplicationwithtilde})} \\
%& = & a \wedge b & \mbox{by definition of }  \wedge.  
%\end{array}
%$$
%Then $\mathbf A$ is meet-semidistributive.
%}

Next, we desire to check that $\mathbf A$ satisfies {\rm($M_2$)}.   %(SD$_\land$ (right)). 
Let us assume then that $\mathbf A \models x \wedge y \approx z \wedge y$, and let $a,b,c \in A$. Hence,
%$$
%\begin{array}{lcll}
%(a \vee c) \wedge b & = & ((a' \to c'')'' \to b')' & \mbox{by definition of } \vee \mbox{ and } \wedge  \\
%& = & ((a' \to c) \to b')' & \mbox{by (I$_{2,0})$} \\
%& = & ((b'' \to a') \to (c \to b')')'' & \mbox{by (I)} \\
%& = & ((b'' \to a') \to (c \wedge b))'' &  \mbox{by definition of }  \wedge  \\
%& = & ((b'' \to a') \to (a \wedge b))'' & \mbox{by hypothesis} \\
%& = & ((b'' \to a') \to (a \to b')')''  & \mbox{by definition of }  \wedge  \\
%& = & ((a' \to a) \to b')' & \mbox{by (I)} \\
%& = & (a \to b')' & \mbox{by Lemma \ref{general_properties_equiv} (\ref{LeftImplicationwithtilde})} \\
%& = & a \wedge b & \mbox{by definition of }  \wedge.  
%\end{array}
%$$
\noindent $(a \vee c) \wedge b $
$\overset{  \ {\rm def\  of }\  \vee  \ {\rm and }\  \wedge  
}{=}  ((a' \to c'')'' \to b')' $
$\overset{  (I_{2,0})
}{=}  ((a' \to c) \to b')' $
$\overset{  (I) 
}{=}  ((b'' \to a') \to (c \to b')')'' $
$\overset{   \ {\rm def\  of }\   \wedge  
}{=}  ((b'' \to a') \to (c \wedge b))'' $
$\overset{ hyp 
}{=}  ((b'' \to a') \to (a \wedge b))'' $
$\overset{  \ {\rm def\  of }\  \wedge  
}{=}  ((b'' \to a') \to (a \to b')')''  $
$\overset{  (I) 
}{=}  ((a' \to a) \to b')' $
$\overset{   \ref{general_properties_equiv} (\ref{LeftImplicationwithtilde}) 
}{=}  (a \to b')' $
$\overset{  \ {\rm def\  of }\   \wedge.  
}{=}  a \wedge b $. 
Then $\mathbf A$ is meet-semidistributive.
 
\end{Proof}

\begin{Theorem} \label{Theorem_SDinI20}
Let $\mathbf A \in \mathcal I_{2,0}$. Then $\mathbf A_{mj}$ is semidistributive.
\end{Theorem}	
	
\begin{Proof}
%{\color{red}By Theorem \ref{Theorem_070220_01}, $\mathbf A_{mj}$ is meet-semi-distributive. It is easy to see that  {\rm($J_1$)} 
%and {\rm($J_2$)} 
%follow from (SD$_\land$ (left)) and (SD$_\land$ (right)) \cite[Theorem 7.1]{cornejo2015implication}.}
%
%JuanMa, the blue text below is a slight modification of the red text given above.  Check the blue and if you agree with it, then delete the red text. 

By Theorem \ref{Theorem_070220_01}, $\mathbf A_{mj}$ is meet-semi-distributive. It is easy to see that {\rm($J_1$)} 
and {\rm($J_2$)} 
follow from {\rm($M_1$)} 
and {\rm($M_2$)} in view of  \cite[Theorem 7.1]{cornejo2015implication}. Hence $\mathbf A$ is semidistributive.

\end{Proof}

\section{Semidistributivity of I-zroupoids}

%{\color{red}Sanka, the proof of Theorem 6.1 has a very long proof. It reminds me the associative one. So My proposal is to left this result for another paper. What do you think? }

\indent In this section, we prove the main theorem of this paper.  For this, we need one more crucial result proved in \cite{cornejo2016derived}.

\begin{Theorem}[Transfer Theorem] {\rm \cite{cornejo2016derived}}  \label{Theo_identities_I20_to_I} \label{Transfer_Theorem}
	
	Let $t_i(\overline x), i= 1, \cdots, 6$.
	be terms, where $\overline x$ denotes the sequence $\langle x_1, \cdots x_n \rangle$, $x_i$ being varaibles.  Let $\mathcal V$ be a subvariety of $\mathcal I$.  \\ 
	If 
	$$\mathcal V  \cap \mathcal I_{2,0} \models (t_1(\overline x) \to t_2(\overline x)) \to t_3(\overline x) \approx (t_4(\overline x) \to t_5(\overline x)) \to t_6(\overline x),$$ 
	then 
	$$\mathcal V\models (t_1(\overline x) \to t_2(\overline x)) \to t_3(\overline x) \approx (t_4(\overline x) \to t_5(\overline x)) \to t_6(\overline x).$$
\end{Theorem}

We are now ready to present our first main result of this paper.  It was quite surprising to the authors that this result holds.

\begin{Theorem} \label{Main}
Let $\mathbf A \in \mathcal I$. Then $\mathbf A_{mj}$ is semidistributive.
\end{Theorem}

\begin{Proof}
	Apply Theorem \ref{Theorem_SDinI20} and Theorem \ref{Theo_identities_I20_to_I}.	
\end{Proof}

\, \,
The following corollaries are immediate.

\begin{Corollary}
If $\mathbf A$ is a symmetric I-zroupoid, then $\mathbf A_{mj}$ is semidistributive.
\end{Corollary}

%{\color{red} I think that we need to recall or put the definition of implication semigroup here.}

An I-zroupoid $\mathbf A$ is called an implication semigroup if $\mathbf A$ satisfies the associative identity:\\
$ x \to (y \to z) \approx (x \to y) \to z$.\\
A complete description of the lattice of subvarieties of the variety of implication semigroups is given in \cite{Implsem}.

\begin{Corollary}
If $\mathbf A$ is an implication semigroup, then $\mathbf A_{mj}$ is semidistributive.
\end{Corollary}

%{\color{red} I think that we need to recall or put the definition of Birkhoff bisemigroup here.}

In \cite{cornejo2020birkhoff}, we generalized the notion of Birkhoff systems to Birkhoff bisemigroups.
Recall from \cite{cornejo2020birkhoff} that a bisemigroup $\mathbf A$ is a Birkhoff bisemigroup if $\mathbf A$ satisfies the Birkhof identity:\\

 (BR)  $x \land (x \lor y) \approx x \lor(x \land y)$.\\

The following corollary is immediate from Theorem \ref{Main} and one of  main results of \cite{cornejo2020birkhoff}. 
 
\begin{Corollary} \label{CorA}
If $\mathbf A \in \mathcal I$, then $\mathbf A_{mj}$ is a semidistributive, Birkhoff bisemigroup.
\end{Corollary}

\section{Whitman Property}

\

It is a well-known result in lattice theory, proved by P. Whitman \cite{Wh41} that every free lattice satisfies the following property (W): 
\[ x \land y \leq z \lor u \Rightarrow x \leq z \lor u$ \  or \ $ y < z \lor u $ \ or \ $x \land y \leq z$ \ or \ $ x \land y \leq u. \] 

(W) is now known as Whitman Property.
%{\color{red}A lattice L satisfies Whitman's property whenever\\
%(W) for all $a, b, c \in L$, if $a \land b <  c \lor d$, then  $a < c \lor d$ \ or $b < c \lor d$ \ or $a \land b < c$ \  or $a \land b < d$. \\
%}

As a matter of fact, (W) is one of four conditions used by P. M. Whitman \cite{Wh41} to characterize free lattices.  The Whitman condition was studied, for example, in \cite{freese1995free} and \cite{day1977splitting}.   Also,
B. Jonsson and J. E. Kiefer \cite{JoKi62} investigated finite lattices which satisfy
semi-distributivity and  the Whitman Property (W).
%(W) $ x \land y \leq z \lor u \Rightarrow x \leq z \lor u$ and  $x \land y \leq z$ or  $x \land y \leq u$. (W).
%It is a well-known result due to P. Whitman [4] that every free lattice satisfies (W).'' \\
Since the variety $\mathcal I$ satisfies semidistributivity, it was only natural to wonder if $\mathcal I$ or any of its subvarieties satisfies (W).

 The following algebra shows that (W) fails in $\mathcal I$ (for x=2, y=3, z=1, u=4):

\medskip

\begin{tabular}{r|rrrrrrr}
	$\to$: & 0 & 1 & 2 & 3 & 4 & 5 & 6\\
	\hline
	0 & 5 & 5 & 5 & 5 & 5 & 5 & 5 \\
	1 & 2 & 2 & 2 & 5 & 2 & 5 & 2 \\
	2 & 1 & 1 & 2 & 3 & 6 & 5 & 6 \\
	3 & 4 & 6 & 2 & 3 & 4 & 5 & 6 \\
	4 & 3 & 3 & 5 & 3 & 3 & 5 & 3 \\
	5 & 0 & 1 & 2 & 3 & 4 & 5 & 6 \\
	6 & 6 & 6 & 2 & 3 & 6 & 5 & 6
\end{tabular}

\ \\

We were surprised to find that, indeed, there are subvarieties of $\mathcal I$ 
that do satisfy (W).\\

	Let $\mathcal{MEJ}$ (``meet equals join'') denote the subvariety of $\mathcal I$ defined by the identity:\\
	$ x \land y \approx x \lor y. $  %The relationships of $\mathcal{MEJ}$ with other known subvarieties of $\mathcal{I}$ are investigated in \cite{cornejo2020varieties}. 
	 Some properties of  $\mathcal{MEJ}$ and its relationships with other known subvarieties of $\mathcal{I}$ are investigated in \cite{cornejo2020varieties}.  We just mention here that the variety of implication semigroups mentioned at the end of the previous section is, in fact, a subvariety of $\mathcal{MEJ}$, as proved in \cite{cornejo2020varieties}. 
	\\
	
%	{\color{red} We may mention here why is important the subvariety $\mathcal{MEJ}$, like a little talk.}\\
%	
%	{\color{blue}JuanMa,  I am not sure what to say here.  If you have some idea what we should put here, then, please type it here.  I will also continue thinking about this.}
	
%Recall $ x \leq_{\land} y$ iff x $\land y = x. $

\begin{lema} \label{MEJ_properties}
	Let $\mathbf A \in \mathcal{MEJ}$. Then $\mathbf A$ satisfies:
	\begin{enumerate}[{\rm (1)}]
		\item $(x \to y')' \approx (x' \to y'')''$, \label{090620_01}
		\item $0 \to [(x \to y) \to z] \approx (x \to y) \to z,$ \label{090620_02}
		\item $0 \approx 0',$ \label{090620_03}
		\item $0 \to [x \to (y \to z)] \approx (x \to y')' \to z,$ \label{090620_04}
		%\item $[0 \to \{x \to (y \to z)\}] \wedge z \approx 0 \to [x \to (y \to z)].$ \label{090620_05}
		\item $[0 \to \{x \to (y \to z)\}] \leq z.$ %\approx 0 \to [x \to (y \to z)].$ 
		\label{090620_05} 
	\end{enumerate}
\end{lema}

\begin{Proof}  Let $a,b,c \in A$. Then
	\begin{itemize}
		\item[(\ref{090620_01})] 
%		 $$
%		\begin{array}{lcll}
%		(a \to b')' & = & a \wedge b & \mbox{by definition of } \land \\
%		& = & a \vee b & \mbox{since }  x \land y \approx x \lor y \\
%		& = & (a' \wedge b')' & \mbox{by definition of } \lor \\
%		& = & (a' \to b'')'' & \mbox{by definition of } \land. 
%		\end{array}
%		$$
\noindent $		(a \to b')' $
$\overset{  def\ of\ \land 
}{=}  a \wedge b $
$\overset{   x \land y \approx x \lor y 
}{=}  a \vee b $
$\overset{  def\ of\  \lor 
}{=}  (a' \wedge b')' $
$\overset{  def\ of\  \land 
}{=}  (a' \to b'')'' $.
		
		\item[(\ref{090620_02})] 
%		 $$
%		\begin{array}{lcll}
%		0 \to [(a \to b) \to c] & = & 0 \to [(c' \to a) \to (b \to c)']' & \mbox{by (I)} \\
%		& = & 0'' \to [(c' \to a) \to (b \to c)']' & \mbox{since } 0 \approx 0'' \\
%		& = & [0'' \to ((c' \to a) \to (b \to c)')']'' & \mbox{by Lemma \ref{Lemma_300315_01} (\ref{060415_01})} \\
%		& = & [0''' \to ((c' \to a) \to (b \to c)')'']''' & \mbox{by (\ref{090620_01})} \\
%		& = & [0' \to ((c' \to a) \to (b \to c)')'']''' & \mbox{since } 0 \approx 0'' \\
%		& = &  [((c' \to a) \to (b \to c)')'']''' &  \mbox{by Lemma \ref{general_properties_equiv} (\ref{TXX}) and Theorem \ref{Transfer_Theorem}} \\
%		& = & ((c' \to a) \to (b \to c)')' & \mbox{by Lemma \ref{Lemma_300315_01} (\ref{060415_01})} \\
%		& = & (a \to b) \to c & \mbox{by (I)} 
%		\end{array}
%		$$
\noindent $		0 \to [(a \to b) \to c] $
$\overset{  (I) 
}{=}  0 \to [(c' \to a) \to (b \to c)']' $
$\overset{  0 \approx 0'' 
}{=}  0'' \to [(c' \to a) \to (b \to c)']' $
$\overset{   \ref{Lemma_300315_01} (\ref{060415_01}) 
}{=}  [0'' \to ((c' \to a) \to (b \to c)')']'' $
$\overset{  (\ref{090620_01}) 
}{=}  [0''' \to ((c' \to a) \to (b \to c)')'']''' $
$\overset{  0 \approx 0'' 
}{=}  [0' \to ((c' \to a) \to (b \to c)')'']''' $
$\overset{    \ref{general_properties_equiv} (\ref{TXX}) and  \ref{Transfer_Theorem} 
}{=}   [((c' \to a) \to (b \to c)')'']''' $
$\overset{   \ref{Lemma_300315_01} (\ref{060415_01}) 
}{=}  ((c' \to a) \to (b \to c)')' $
$\overset{  (I) 
}{=}  (a \to b) \to c $.

		\item[(\ref{090620_03})] 
%		 $$
%		\begin{array}{lcll}
%	0'	& = & 0 \to 0 & \mbox{} \\
%		& = & 0 \to 0'' & \mbox{since } 0 \approx 0'' \\
%		& = & 0 \to (0' \to 0) & \mbox{} \\
%		& = & 0' \to 0 & \mbox{by (\ref{090620_02})} \\
%		& = & 0'' & \mbox{} \\
%		& = & 0 & \mbox{since } 0 \approx 0''
%		\end{array}
%		$$
\noindent $	0'	$
$\overset{  
}{=}  0 \to 0 $
$\overset{  0 \approx 0'' 
}{=}  0 \to 0'' $
$\overset{  
}{=}  0 \to (0' \to 0) $
$\overset{  (\ref{090620_02}) 
}{=}  0' \to 0 $
$\overset{  
}{=}  0'' $
$\overset{  0 \approx 0''
}{=}  0 $.
		
		\item[(\ref{090620_04})] 
%	 $$
%	\begin{array}{lcll}
%0 \to (a \to (b \to c))	& = & 0' \to (a \to (b \to c)) & \mbox{by (\ref{090620_03})} \\
%	& = & 0' \to (a'' \to (b \to c))  & \mbox{by Theorem \ref{Transfer_Theorem}} \\
%	& = & 0 \to (a'' \to (b \to c)) & \mbox{by (\ref{090620_03})} \\
%	& = & 0 \to ((a' \to 0) \to (b \to c)) & \mbox{ } \\%by (\ref{090620_03})} \\
%	& = & a'' \to (b \to c) & \mbox{by (\ref{090620_02})} \\
%	& = & (a'' \to b')' \to c & \mbox{by Lemma \ref{general_properties} (\ref{140715_20}) and Theorem \ref{Transfer_Theorem}} \\
%	& = & (a \to b')' \to c & \mbox{by Theorem \ref{Transfer_Theorem}} 
%	\end{array}
%	$$
\noindent $0 \to (a \to (b \to c))	$
$\overset{  (\ref{090620_03}) 
}{=}  0' \to (a \to (b \to c)) $
$\overset{   \ref{Transfer_Theorem} 
}{=}  0' \to (a'' \to (b \to c))  $
$\overset{  (\ref{090620_03}) 
}{=}  0 \to (a'' \to (b \to c)) $
$\overset{  
}{=}  0 \to ((a' \to 0) \to (b \to c)) $
$\overset{  (\ref{090620_02}) 
}{=}  a'' \to (b \to c) $
$\overset{   \ref{general_properties} (\ref{140715_20}) and  \ref{Transfer_Theorem} 
}{=}  (a'' \to b')' \to c $
$\overset{   \ref{Transfer_Theorem} 
}{=}  (a \to b')' \to c $.
	
	\item[(\ref{090620_05})] 
% $$
%\begin{array}{lcll}
%[0 \to (a \to (b \to c))] \wedge c & = & [[0 \to (a \to (b \to c))] \to c']' & \mbox{by definition of } \land \\
%& = & [[(a \to b')' \to c] \to c']' & \mbox{by (\ref{090620_04})} \\
%& = & [[(a \to b')' \to c''] \to c']' & \mbox{by Theorem \ref{Transfer_Theorem}} \\
%& = & [[(a \to b')' \to 0'] \to c']' & \mbox{by Lemma \ref{general_properties2} (\ref{281014_05})} \\
%& = & [(a \to b')'' \to c']' & \mbox{by (\ref{090620_03})} \\
%& = & [(a \to b')''' \to c'']'' & \mbox{by (\ref{090620_01})} \\
%& = & (a \to b')' \to c & \mbox{by Theorem \ref{Transfer_Theorem}} \\
%& = & 0 \to (a \to (b \to c)) & \mbox{by (\ref{090620_04})}
%\end{array}
%$$
\noindent $[0 \to (a \to (b \to c))] \wedge c $
$\overset{  def\  of\   \land 
}{=}  [[0 \to (a \to (b \to c))] \to c']' $
$\overset{  (\ref{090620_04}) 
}{=}  [[(a \to b')' \to c] \to c']' $
$\overset{   \ref{Transfer_Theorem} 
}{=}  [[(a \to b')' \to c''] \to c']' $
$\overset{   \ref{general_properties2} (\ref{281014_05}) 
}{=}  [[(a \to b')' \to 0'] \to c']' $
$\overset{  (\ref{090620_03}) 
}{=}  [(a \to b')'' \to c']' $
$\overset{  (\ref{090620_01}) 
}{=}  [(a \to b')''' \to c'']'' $
$\overset{   \ref{Transfer_Theorem} 
}{=}  (a \to b')' \to c $
$\overset{  (\ref{090620_04})
}{=}  0 \to (a \to (b \to c)) $.
		
	\end{itemize}	
\end{Proof}

\begin{lema} \label{MEJ_properties2}
	Let $\mathbf A \in \mathcal{MEJ}$. Let $a,b,c,d \in  A$ such that $[(a \to b')' \to (c \to d')'']' = (a \to b')'$. Then:
	\begin{enumerate}[{\rm (1)}]
		\item $(a \to b') \to (c' \to d) = a' \to b$ \label{090620_06}
		\item $0 \to (a' \to b) \leq d.$ \label{090620_07}
	\end{enumerate}
\end{lema}

\begin{Proof}
	\begin{itemize}
	\item[(\ref{090620_06})] 
%	 $$
%	\begin{array}{lcll}
%	(a \to b') \to (c' \to d) & = & (a \to b') \to (c' \to d)'' & \mbox{by Lemma \ref{Lemma_300315_01} (\ref{060415_01})} \\
%	& = & (a \to b') \to (c'' \to d')''' & \mbox{by Lemma \ref{MEJ_properties} (\ref{090620_01})} \\
%	& = & (a \to b') \to (c \to d')' & \mbox{by Theorem \ref{Transfer_Theorem}} \\
%	& = & [(a \to b') \to (c \to d')']'' & \mbox{by Lemma \ref{Lemma_300315_01} (\ref{060415_01})} \\
%	& = & [(a \to b')' \to (c \to d')'']''' & \mbox{by Lemma \ref{MEJ_properties} (\ref{090620_01})} \\
%	& = & [(a \to b')' \to (c \to d')'']' & \mbox{by Theorem \ref{Transfer_Theorem}} \\
%	& = & (a \to b')' & \mbox{by hypothesis} \\
%	& = & (a' \to b'')'' & \mbox{by Lemma \ref{MEJ_properties} (\ref{090620_01})} \\
%	& = & a' \to b & \mbox{by Theorem \ref{Transfer_Theorem}} 
%	\end{array}
%	$$
\noindent $	(a \to b') \to (c' \to d) $
$\overset{   \ref{Lemma_300315_01} (\ref{060415_01}) 
}{=}  (a \to b') \to (c' \to d)'' $
$\overset{   \ref{MEJ_properties} (\ref{090620_01}) 
}{=}  (a \to b') \to (c'' \to d')''' $
$\overset{   \ref{Transfer_Theorem} 
}{=}  (a \to b') \to (c \to d')' $
$\overset{   \ref{Lemma_300315_01} (\ref{060415_01}) 
}{=}  [(a \to b') \to (c \to d')']'' $
$\overset{   \ref{MEJ_properties} (\ref{090620_01}) 
}{=}  [(a \to b')' \to (c \to d')'']''' $
$\overset{   \ref{Transfer_Theorem} 
}{=}  [(a \to b')' \to (c \to d')'']' $
$\overset{ hyp 
}{=}  (a \to b')' $
$\overset{   \ref{MEJ_properties} (\ref{090620_01}) 
}{=}  (a' \to b'')'' $
$\overset{   \ref{Transfer_Theorem} 
}{=}  a' \to b $.
	
	\item[(\ref{090620_07})] 
%	 $$
%	\begin{array}{lcll}
%0 \to (a' \to b)	& = & 0 \to ((a \to b') \to (c' \to d)) & \mbox{by (\ref{090620_06})} \\
%	& \leq & d & \mbox{by Lemma \ref{MEJ_properties} (\ref{090620_05}) with } x:=  a \to b',\\
%	&  &  & y:= c', z:= d. 
%	\end{array}
%	$$
\noindent $0 \to (a' \to b)	$
$\overset{  (\ref{090620_06}) 
}{=}  0 \to ((a \to b') \to (c' \to d)) $
$\overset{}{\leq}  d $ by Lemma \ref{MEJ_properties} (\ref{090620_05}) with  $x:=  a \to b', y:= c', z:= d.$
%$\overset{ y:= c', z:= d. 
%}{=}   $
	\end{itemize}
\end{Proof}

We now present our second main result of the paper.

\begin{Theorem} \label{SecondMain}
The variety $\mathcal{MEJ}$ satisfies the Whitman Property (W). 
\end{Theorem}

\begin{Proof}
Let $\mathbf A \in \mathcal{MEJ}$ and $a,b,c,d \in A$ such that $a \wedge b \leq c \vee d$. Since the identity $x \wedge y \approx x \vee y$ holds in $\mathcal{MEJ}$ we have that $a \wedge b \leq c \wedge d$. Then
%$$
%\begin{array}{lcll}
%[(a \to b')' \to (c \to d')'']' %& = & (a \to b')' \wedge (c \to d')' & \mbox{by definition of } \land \\
%& = & (a \wedge b) \wedge (c \wedge d) & \mbox{by definition of } \land \\
%& = & a \wedge b & \mbox{since } a \wedge b \leq c \wedge d \\
%& = & (a \to b')' & \mbox{by definition of } \land. 
%\end{array}
%$$
\noindent $[(a \to b')' \to (c \to d')'']'$ %$
%$\overset{  def\  of\   \land 
%}{=}  (a \to b')' \wedge (c \to d')' $
$\overset{  def\  of\   \land 
}{=}  (a \wedge b) \wedge (c \wedge d) $
$\overset{  a \wedge b \leq c \wedge d 
}{=}  a \wedge b $
$\overset{  def\  of\   \land 
}{=}  (a \to b')' $. 
Hence, by Lemma \ref{MEJ_properties2} (\ref{090620_07}), 
\begin{equation} \label{equation_090620_01}
0 \to (a' \to b) \leq d.
\end{equation}
Therefore,
%$$
%\begin{array}{lcll}
%a \wedge b & = & (a \to b')' & \mbox{by definition of } \land \\
%& = & (a' \to b'')'' & \mbox{by Lemma \ref{MEJ_properties} (\ref{090620_01})} \\
%& = & a' \to b & \mbox{by Theorem \ref{Transfer_Theorem}} \\
%& = & 0 \to (a' \to b) & \mbox{by Lemma \ref{MEJ_properties} (\ref{090620_02})} \\
%& \leq & d & \mbox{by (\ref{equation_090620_01})} 
%\end{array}
%$$
\noindent $a \wedge b $
$\overset{  def\  of\  \land 
}{=}  (a \to b')' $
$\overset{   \ref{MEJ_properties} (\ref{090620_01}) 
}{=}  (a' \to b'')'' $
$\overset{   \ref{Transfer_Theorem} 
}{=}  a' \to b $
$\overset{   \ref{MEJ_properties} (\ref{090620_02}) 
}{=}  0 \to (a' \to b) $
$\overset{  (\ref{equation_090620_01}) 
}{\leq}  d $.
\end{Proof}

It should be remarked here that the relation $\leq$ is not, in general, a partial order on algebras in $\mathcal I$.  In fact, it was shown in \cite{cornejo2016order} that the variety $\mathcal{I}_{2, 0}$ is the maximal subvariety of 
$\mathcal I$ with respect to the property that $\leq$ is a partial order in it.  Thus we can state the following corollary where the relation $\leq$ is indeed a partial order. 
 
The following corollary is immediate from Corollary \ref{CorA} and Theorem \ref{SecondMain}. 
\begin{Corollary}
If $\mathbf A \in \mathcal{MEJ} \cap \mathcal{I}_{2, 0}$, then $\mathbf A_{mj}$ is a semidistributive, Birkhoff bisemigroup which satisfies Whitman Property (W).
\end{Corollary}

\section{Concluding Remarks}

It is well-known that the variety of semilattices is congruence-semi-distributive. 
A proof is given in \cite{Pa64}.
%the paper: Papert, Dona. Congruence relations in semi-lattices. J. London Math. Soc. 39 1964 723-729.   
It is also well-known that the variety of De Morgan algebras is congruence-distributive.  Since the variety of implication zroupoids contains both the variety of $\lor$-semilattices with 0 and the variety of De Morgan algebras, the following problem arises naturally.\\

{\bf PROBLEM 1}:  Is the variety of implication zroupoids congruence-semidistributive? \\

Although the classes of meet semi-distibutive bisemigroups, join semi-distibutivie bisemigroups,  and semi-distibutive bisemigroups are natural extensions of the corresponding classes of lattices,  they do not, to the best of our knowledge, seem to have been investigated in the literature so far.  The main result of this paper certainly seem to warrant such an investigation.  To facilitate such a study, we mention the following open problem which naturally arises.\\

{\bf PROBLEM 2}  Investigate the quasivarietties of meet semi-distibutive bisemigroups, join semi-distibutivie bisemigroups and semi-distibutive bisemigroups.  \\

Let $\mathbf A =\langle A, \lor, \land \rangle$ be a bisemigroup.  $\mathbf A$ is $\land$-idempotent if $\mathbf A$ satisfies the 
$\land$-idempotent identity: $x \land x \approx x$. 
$\lor$-idempotent bisemigroups are defined dually.  A bisemigroup $\mathbf A$ is idempotent if it is both 
$\land$-idempotent and $\lor$-idempotent.
$\mathbf A$ is left [right] mj-distributive if $\mathbf A \models x \land (y \lor z) \approx (x \land y) \lor (x \land z)$ [$(y \lor z) \land x \approx (y \land x) \lor (z \land x)$].  Left [right] jm-distributive bisemigroups are defined dually.  $\mathbf A$ is mj-distributive if it is both left mj-distributive and right mj-distributive.
jm-distributive bisemigroups are defined dually.  $\mathbf A$ is distributive if it is both mj-distributive and jm-distributive.\\

We conclude this paper with the following (easy) observation.

\begin{Theorem}
Let $\mathbf A$ be a left mj-distributive, $\lor$-idempotent bisemigroup. 
%\begin{itemize}
%	\item $x \land (y \lor z) \approx (x \land y) \lor (x \land z)$ and
%	\item $x \lor x \approx x$.
%\end{itemize}
Then
%$x \land  y = x \land  z \Rightarrow  x \land  y = x \land (y \lor z)$.
$\mathbf A$ is meet-semidistributive.
In particular, every distributive, idempotent bisemigroup is semi-distributive.
\end{Theorem} 

%\begin{comment}
%1 x ^ y = x ^ z -> x ^ y = x ^ (y v z) # label(non_clause) # label(goal).  [goal].
%2 x ^ (y v z) = (x ^ y) v (x ^ z).  [assumption].
%3 (x ^ y) v (x ^ z) = x ^ (y v z).  [copy(2),flip(a)].
%4 x v x = x.  [assumption].
%5 c1 ^ c3 = c1 ^ c2.  [deny(1)].
%6 c1 ^ (c2 v c3) != c1 ^ c2.  [deny(1)].
%7 (c1 ^ x) v (c1 ^ c2) = c1 ^ (x v c3).  [para(5(a,1),3(a,1,2))].
%8 c1 ^ (x v c2) = c1 ^ (x v c3).  [para(3(a,1),7(a,1))].
%9 c1 ^ (x v c3) = c1 ^ (x v c2).  [copy(8),flip(a)].
%10 c1 ^ (c2 v c2) != c1 ^ c2.  [para(9(a,1),6(a,1))].
%11 c1 ^ c2 != c1 ^ c2.  [para(4(a,1),10(a,1,2))].
%12 $F.  [copy(11),xx(a)].

%============================== end of proof ==========================
%\end{comment}

%We conclude this paper with the following open problem. \\

%{\bf PROBLEM:} Investigate semi-distributive Birkhoff bisemigroups. %in particular, describe the structure of the lattice of subvarieties of the variety of Birkhoff bisemigroups.

%\begin{Theorem}
%Let $\mathbf{A} \in \mathcal I$.  Then \\
% $\mathbf{A_{mj}}$ is semidistributive.\
%\end{Theorem}
%
%\begin{Proof}
%It suffices to prove that (MSD) holds in $\mathbf{A_{mj}}$, since (JSD) follows from (MSD).
%\end{Proof}

\section*{Acknowledgements} 

 The first author wants to thank the institutional support of CONICET (Consejo Nacional de Investigaciones Cient\'ificas y T\'ecnicas) and Universidad Nacional del Sur.  The authors also wish to acknowledge that \cite{Mc} was a useful tool during the research phase of this paper.\\

\noindent {\bf Compliance with Ethical Standards:}\\ 

{\bf Conflict of Interest}: The first author declares that he has no conflict of interest. The second author declares that he has no conflict of interest.\\

{\bf Ethical approval}: This article does not contain any studies with human participants or animals performed by any of the authors.\\

 {\bf Funding:}  
	The work of Juan M. Cornejo was supported by CONICET (Consejo Nacional de Investigaciones Cientificas y Tecnicas) and Universidad Nacional del Sur.

\newpage
\vskip 1.5cm

\noindent {\sc Juan M. Cornejo}\\
Departamento de Matem\'atica\\
Universidad Nacional del Sur\\
Alem 1253, Bah\'ia Blanca, Argentina\\
INMABB - CONICET\\

\noindent jmcornejo@uns.edu.ar

\vskip 1.4cm

\noindent {\sc Hanamantagouda P. Sankappanavar}\\
Department of Mathematics\\
State University of New York\\
New Paltz, New York 12561\\
U.S.A.\\

\noindent sankapph@newpaltz.edu

\end{document}